\documentclass[12pt]{amsart}
\usepackage{amsmath,amssymb,latexsym,amsfonts,amsthm,amscd}%
\usepackage{enumerate}
\usepackage[mathscr]{eucal}

\numberwithin{equation}{section}

\theoremstyle{plain}
\newtheorem{thm}{Theorem}[section]
\newtheorem{prop}[thm]{Proposition}
\newtheorem{lem}[thm]{Lemma}
\newtheorem{coro}[thm]{Corollary}

\theoremstyle{definition}
\newtheorem{defi}[thm]{Definition}

\newtheorem{assu}[thm]{Assumption}

\newtheorem{rem}[thm]{Remark}
\newtheorem{exa}[thm]{Example}

\def\E{{\mathscr{E}}}
\def\F{{\mathscr{F}}}
\def\R{ {\mathbf R} }
\def\C{ {\mathbf C} }
\def\Z{ {\mathbf Z} }
\def\N{ {\mathbf N} }
\def\T{ {\mathbf T} }

\def\sP{ \mathscr{P} }

\def\bP{ \boldsymbol{P} }
\def\bQ{ \boldsymbol{Q} }
\def\bE{ \boldsymbol{E} }

\def\sG{ \mathscr{G} }
\def\sH{ \mathscr{H} }

\def\fS{ \mathfrak{S} }

\def\ext{ {\rm ex} (\sP_\mu) }
\def\sgn{ {\rm sgn} }
\def\eps{ \varepsilon }

\renewcommand{\hat}{\widehat}
\renewcommand{\tilde}{\widetilde}

\newcommand{\unit}{e} 

\newcommand{\cbra}[1]{\left( #1 \right)}
\newcommand{\kbra}[1]{\left\{ #1 \right\}}
\newcommand{\ebra}[1]{\left[ #1 \right]}

\newcommand{\Cz}{$ ({\bf C0}) $}
\newcommand{\Ca}{$ ({\bf C1}) $}
\newcommand{\Cb}{$ ({\bf C2}) $}
\newcommand{\Cc}{$ ({\bf C3}) $}

\title[Stochastic equation on compact groups]{Stochastic equation on compact groups in discrete negative time}
\author[J.~Akahori]{Jir\^o Akahori}
\address{Jir\^o Akahori, Department of Mathematical Sciences, 
Ritsumeikan University, 
1-1-1 Noji-higashi, Kusatsu, Shiga, 525-8577, Japan}
\email{akahori@se.ritsumei.ac.jp}
\author[C.~Uenishi]{Chihiro Uenishi}
\address{Chihiro Uenishi, Cadem Corporation, Ltd., 
Shin-Yokohama Office, 
Houei-Shin-Yokohama Bldg., 2-14-9, 
Shin-Yokohama, Kouhoku-Ku, Yokohama, 222-0033, Japan}
\author[K.~Yano]{Kouji Yano}
\address{Kouji Yano, 
Research Institute for Mathematical Sciences, 
Kyoto University, 
Kyoto, 606-8502, Japan}
\email{yano@kurims.kyoto-u.ac.jp}
\dedicatory{Dedicated to Professor~Shinzo~Watanabe on the~occasion of his~70th~birthday \\
and to Professor~Yoichiro~Takahashi on the~occasion of his~60th~birthday}
\thanks{This research was supported by 
Open Research Center Project for Private Universities: matching fund 
subsidy from MEXT, 2004-2008.}

\begin{document}

\maketitle
\centerline{March 3, 2006}

\begin{abstract}
In this paper 
a stochastic equation on compact groups in discrete negative time 
is studied. 
This is closely related to Tsirelson's stochastic differential equation, 
of which any solution is non-strong. 
How the group action reflects on the set of solutions is investigated. 
It is applied to generalize Yor's result 
and give a necessary and sufficient condition 
for existence of a strong solution and for uniqueness in law. 
\end{abstract}

\section{Introduction}

In contrast with that of ordinary ones, 
the theory of stochastic differential equations 
has the distinguished notions 
of a {\em strong solution} and two uniqueness properties: 
{\em pathwise uniqueness} and {\em uniqueness in law}. 
Therefore we have the following four cases: 
\begin{equation}
\text{
\begin{tabular}{c||c|c}
                                 & unique in law   & non-unique in law \\
\hline
\hline
$ \exists $ strong solution      & \Cz             & \Cb \\
\hline
$ \not \exists $ strong solution & \Ca             & \Cc 
\end{tabular}
}
\label{table}
\end{equation}
The celebrated theorem of T.~Yamada and S.~Watanabe (\cite{YW}) 
is stated as follows: 
{\em Pathwise uniqueness implies uniqueness in law 
and then any solution is strong}. 
Except for the trivial cases where there is no solution, 
we may say that pathwise uniqueness implies that the case {\Cz} occurs. 
In most cases the converse is also true (See A.~K.~Zvonkin and N.~V.~Krylov \cite{ZK}). 

B.~Tsirelson (\cite{Tsi}) has presented 
a remarkable example of a stochastic differential equation \eqref{Tsirel SDE} stated below, 
which enjoys uniqueness in law property but has no strong solution: 
In short, the case {\Ca} occurs. 
His equation has deeply been investigated from various viewpoints 
by many researchers: See, for example, \cite{ES}, \cite{Kall} and \cite{LGY}. 

To prove the non-existence of strong solutions, 
Tsirelson introduced a stochastic equation on the torus $ \T = \R / \Z $ 
indexed by discrete negative time: 
\begin{equation}
\eta_k = \eta_{k-1} + \xi_k , 
\quad k \in -\N, 
\label{Tsirel eq}
\end{equation}
where $ \xi = ( \xi_k , k \in -\N ) $ is a Gaussian driving noise 
and $ \eta = ( \eta_k , k \in -\N ) $ is an unknown process. 
We can find its simpler proof, mainly due to N.~V.~Krylov, 
in the literature: 
For example, see \cite[pp.150--151]{LS}, \cite[pp.149--151]{SY} and \cite[pp.195--197]{IW}. 
In these contexts, 
the authors discussed, instead of \eqref{Tsirel eq}, 
a modified equation on the real line $ \R $ as follows: 
\begin{eqnarray}
\eta_k = \alpha (\eta_{k-1})  + \xi_k , 
\quad k \in -\N. 
\label{Yor eq}
\end{eqnarray}
Here $ \alpha (x) $ denotes the fractional part of $ x \in \R $. 

M.~Yor (\cite{Yor}) has studied the stochastic equation \eqref{Yor eq} 
for an {\em arbitrarily} given noise law, which is not necessarily Gaussian. 
In this case we always have a non-strong solution, 
so that the case {\Cz} never occurs. 
He successfully obtained striking results 
to give a necessary and sufficient condition for the trichotomy {\Ca}-{\Cc} 
in terms of the noise law. 
One of his results may roughly be stated as follows. 
To any given noise law $ \mu $ 
there corresponds a subgroup $ \Z_\mu $ of $ \Z $ 
such that the following holds. 
\begin{thm}[Yor \cite{Yor}] \label{Intro YorTh}
{\rm (i)} 
The case {\Ca} occurs iff $ \Z_\mu = \kbra{0} $. 
\\
{\rm (ii)} 
The case {\Cb} occurs iff $ \Z_\mu = \Z $. 
If it is true, then, among the solutions, 
all extremal points are strong and 
the others non-strong. 
\\
{\rm (iii)} 
The case {\Cc} occurs iff $ \kbra{0} \subsetneq \Z_\mu \subsetneq \Z $. 
\end{thm}
This result sounds somehow paradoxical; at least, they are mysterious: 
{\em Uniqueness in law implies that the solution is non-strong, 
and existence of strong solution(s) 
implies that it is non-unique but extremal}. 

\

The purpose of the present paper is 
to generalize Yor's results. 
We study the following stochastic equation on a compact group $ G $ in discrete negative time: 
\begin{equation}
\eta_k = \eta_{k-1} \cdot \xi_k , 
\quad k \in -\N, 
\label{STYE}
\end{equation}
which we call the {\em simple Tsirelson--Yor equation on $ G $}. 
We shall introduce the notions of a solution, a strong solution and uniqueness in law 
following the theory of stochastic differential equations. 
Then we obtain the following result, 
which generalizes Theorems 2 of \cite{Yor}. 

\begin{thm} \label{Intro exist}
For any given noise law $ \mu $, there exists a solution $ \bP^*_\mu $ 
of the equation \eqref{STYE} such that 
all marginal distributions of $ \bP^*_\mu $ are uniform. 
If $ G \neq \kbra{\unit} $, then the solution $ \bP^*_\mu $ is non-strong. 
\end{thm}

This result says that, if $ G \neq \kbra{\unit} $, 
then the case {\Cz} always fails. 
Theorem \ref{Intro exist} will be restated in Theorem \ref{YorTh} 

We shall naturally define a group action of $ G $ on the set of solutions. 
Our key result is roughly stated as follows. 
\begin{thm} \label{Intro action}
The action restricted on the set of extremal points is transitive. 
\end{thm}
This implies that the set of extremal points is exhausted 
by the $ G $-orbit of an arbitrary extremal point $ \bP^o $ 
and is homeomorphic 
to the homogeneous space $ G/H_\mu(\bP^o) $ 
where $ H_\mu(\bP^o) $ is the isotropic subgroup at the point $ \bP^o $. 
The precise statement of Theorem \ref{Intro action} will be given in Theorem \ref{Gal1}. 

Based on Theorem \ref{Intro action} 
and fully employing the representation theory of compact groups, 
we generalize Theorem \ref{Intro YorTh} 
to give a necessary and sufficient condition 
for the trichotomy {\Ca}-{\Cc} 
in terms of two subgroups $ H_\mu(\bP^o) $ and $ H_\mu^s $. 
Here $ H_\mu^s $ is a closed normal subgroup of $ G $ 
which we will define in \eqref{H mu s}. 

\begin{thm} \label{Intro trichotomy}
{\rm (i)} 
The case {\Ca} occurs iff $ H_\mu(\bP^o) = G $ for some (and hence any) extremal point $ \bP^o $. 
\\
{\rm (ii)} 
The case {\Cb} occurs iff $ H_\mu^s = \kbra{0} $. 
If it is true, then, among the solutions, 
all extremal points are strong and 
the others non-strong. 
\\
{\rm (iii)} 
The case {\Cc} occurs iff $ H_\mu(\bP^o) \subsetneq G $ and $ H_\mu^s \supsetneq \kbra{0} $. 
\end{thm}

Theorem \ref{Intro trichotomy} will be proved in section \ref{sec: trichotomy}. 

The subgroups $ H_\mu^s $ and $ H_\mu(\bP^o) $ for extremal points $ \bP^o $ 
are related as follows. 

\begin{thm} \label{Intro isot group}
For any extremal point $ \bP^o $, the inclusion 
\begin{equation}
H_\mu(\bP^o) \subset H_\mu^s 
\end{equation}
holds. 
If $ G $ is abelian, 
then the equality holds for any extremal point $ \bP^o $. 
\end{thm}

We will restate Theorem \ref{Intro isot group} as Theorem \ref{thm: isot group} 
including further information. 
In the case where $ G = \T = \R/\Z $, 
note that the abelian group $ \Z $ is the dual group of $ G $ 
in the sense of the {\em Pontryagin duality}. 
In this terminology, the subgroup $ \Z_\mu $ of $ \Z $ in Theoren \ref{Intro YorTh} 
is exactly the dual group of $ H_\mu^s $. 
In the non-abelian cases, however, we have a typical example 
given in Example \ref{ex: sym group} 
where $ H_\mu(\bP^o) $ is strictly included in $ H_\mu^s $ for any extremal point $ \bP^o $. 

\

Now we recall Tsirelson's stochastic differential equation (SDE in short). 
Let $ (t_k : \, k \in -\N) $ be a decreasing sequence of the interval $ (0,1] $ such that 
$ t_k \to 0 $ as $ k \to -\infty $. 
Define $ A(t,X) $, which is called {\em Tsirelson's drift}, as 
\begin{equation}
A(t,X) = \sum_{k \in -\N} \left\{ \frac{X_{t_{k-1}} - X_{t_{k-2}}}{t_{k-1} - t_{k-2}} \right\} 
1_{ [t_{k-1},t_{k}) }(t) , \quad 0 < t < 1 . 
\end{equation}
Then Tsirelson's SDE is given by 
\begin{equation}
X_t = B_t + \int_0^t A(s,X) ds , \quad 0 \le t \le 1 
\label{Tsirel SDE}
\end{equation}
where $ (B_t) $ is a one-dimensional Brownian motion. 
If $ ((X_t),(B_t)) $ is a solution of the equation \eqref{Tsirel SDE}, 
then the sequences $ (\eta_k) $ and $ (\xi_k) $ defined by 
\begin{equation}
\eta_k = \frac{X_{t_{k}} - X_{t_{k-1}}}{t_{k} - t_{k-1}} , \quad 
\xi_k = \frac{B_{t_{k}} - B_{t_{k-1}}}{t_k - t_{k-1}} , 
\quad k \in -\N 
\end{equation}
satisfy the equation \eqref{Yor eq}. 
Conversely, we can reconstruct the process $ (X_t) $ 
from the processes $ (\eta_k) $ and $ (B_t) $. 
Hence we may say that all properties of the SDE \eqref{Tsirel SDE} 
can be deduced from those of the equation \eqref{Yor eq}. 

To generalize the stochastic equation \eqref{Yor eq}, 
we introduce a class of stochastic equations in discrete negative time, 
which we call {\em Tsirelson--Yor equations}. 
Let $ S $ be a Polish space and $ G $ a compact group. 
We introduce an operation of $ G $ on the state space $ S $ 
through a measurable map $ \psi:G \times S \to S $ 
and consider a measurable map $ \theta:S \to G $. 
Then our Tsirelson--Yor equation is of the form 
\begin{equation}
\eta_k = \psi \cbra{ \theta(\eta_{k-1}),\xi_k } , 
\quad k \in -\N. 
\label{TYE}
\end{equation}
The key to the generalization is the following 
{\em commutation} condition\footnote{The condition was discovered in \cite{Uen}.}: 
\begin{equation}
\theta \cbra{ \psi( g,s ) } = g \cdot \theta(s), \quad g \in G, \ s \in S. 
\label{UeniC}
\end{equation}
This condition allows us to reduce 
the equation \eqref{T-Yor eq} to the equation 
\begin{equation}
\theta(\eta_k) = \theta (\eta_{k-1}) \cdot \theta(\xi_k) , 
\quad k \in -\N. 
\label{T-Yor eq 2}
\end{equation}
If we write $ \hat{\eta}_k=\theta(\eta_k) $ and $ \hat{\xi}_k=\theta(\xi_k) $, 
then the equation \eqref{T-Yor eq 2} 
is exactly the simple Tsirelson--Yor equation on $ G $. 

We will prove in Proposition \ref{reduction} 
that in order to study the Tsirelson--Yor equation \eqref{T-Yor eq 2} with given laws 
it is sufficient to investigate the simple Tsirelson--Yor equation on $ G $. 

We will show in Example \ref{ex: Yor eq} that Yor's equation \eqref{Yor eq} 
is an example of our Tsirelson--Yor equation. 
Proposition \ref{reduction} says that the equation \eqref{Yor eq} 
is essentially `equivalent' to the equation 
\begin{equation}
\alpha (\eta_k) = \alpha (\eta_{k-1}) + \alpha (\xi_k) \quad \mbox{modulo 1}, 
\quad k \in -\N , 
\label{STYE torus}
\end{equation}
which is nothing but the equation \eqref{Tsirel eq}. 

In addition, 
we will show in Example \ref{ex: Tanaka eq} that 
the following equation is also an example: 
\begin{equation}
\eta_k = \sgn \cbra{\eta_{k-1}} \cdot \xi_k , 
\quad k \in -\N 
\label{Tanaka}
\end{equation}
provided that the noise law has no point mass at $ x=0 $. 
Here $ \sgn (x) = 1 $ if $ x \ge 0 $ and $ \sgn (x) = -1 $ if $ x<0 $. 
This equation has been dealt with 
in \cite[Chapter IX, Exercise 3.18]{R-Y} and \cite[pp. 87--88]{LGY}. 
Proposition \ref{reduction} says that 
the equation \eqref{Tanaka} is equivalent 
to the simple Tsirelson--Yor equation on $ \Z/2\Z $: 
\begin{equation}
\sgn \cbra{\eta_k} = \sgn \cbra{\eta_{k-1}} \cdot \sgn \cbra{\xi_k} , 
\quad k \in -\N . 
\label{STYE Z/2Z}
\end{equation}
This equation will be dealt with in Example \ref{ex: Z/2Z}. 

\ 

The present paper is organized as follows. 
In section \ref{sec: pre}, we give the precise definition 
of the simple Tsirelson--Yor equation on a compact group 
and introduce the notions of a solution, a strong solution and uniqueness in law. 
We also prepare some preliminary facts about the set of solutions. 
Section \ref{sec: exist} is devoted to the proof of Theorem \ref{Intro exist}. 
In section \ref{sec: hom sp}, we give a precise statement and the proof 
of Theorem \ref{Intro action}. 
In section \ref{sec: H mu s}, we introduce a subgroup $ H_\mu^s $ 
which appears in Theorem \ref{Intro trichotomy} 
and restate Theorem \ref{Intro isot group} as Theorem \ref{thm: isot group}. 
In section \ref{sec: example}, we investigate the equation \eqref{STYE Z/2Z} in detail 
and consider another typical example of the simple Tsirelson--Yor equation. 
Section \ref{sec: trichotomy} is devoted to the proofs 
of Theorem \ref{Intro trichotomy} and the propositions 
which are given in section \ref{sec: H mu s}. 
The Tsirelson--Yor equations are defined and discussed in section \ref{sec: TYE}. 

\ 

{\bf Acknowledgments}: 
The authors wish to express sincere thanks 
to Professors Marc Yor, Freddy Delbaen, Shinzo Watanabe, 
and Yoichiro Takahashi 
for stimulating discussions and valuable comments. 
They also thank Gen Mano and Hidehisa Alikawa, who kindly informed them 
of basic theorems of the representation theory. 
The second author expresses his hearty thanks to Professors Hiroki Aoki 
and Yoshiaki Kobayashi for their kind tutorial lectures.

\section{Definitions and preliminary facts} \label{sec: pre}

Let $ G $ be a compact group. 
We consider the following stochastic equation on $ G $ in discrete negative time: 
\begin{equation}
\eta_k = \eta_{k-1} \cdot \xi_k , 
\quad k \in -\N. 
\tag{\ref{STYE}}
\end{equation}

\begin{defi}
The stochastic equation \eqref{STYE} is called 
the {\em simple Tsirelson-Yor equation} on $ G $, 
which will be abbreviated by ``STYE". 
\end{defi}

For an arbitrary sequence $ \mu = (\mu_k , k \in -\N) $ 
of Borel probability measures $ \mu_k $ on $ G $, 
we consider the equation \eqref{STYE} 
with a general noise law $ \mu $. 

Let us give the precise definition of a solution of \eqref{STYE}. 
We denote by $ \eta = ( \eta_k , k \in -\N ) $ 
the coordinate mapping process: 
$ \eta_k(\omega)=\omega(k) $ for $ \omega \in G^{-\N} $ and $ k \in -\N $. 
Set 
\begin{equation}
\xi_k = (\eta_{k-1})^{-1} \cdot \eta_k , 
\quad k \in -\N. 
\label{STYE2}
\end{equation}
Let $ \F^{\eta}_k $ and $ \F^{\xi}_k $ for $ k \in -\N $ 
denote 
\begin{equation}
\F^{\eta}_k := \sigma (\eta_k, \eta_{k-1}, \ldots) , 
\quad k \in -\N 
\end{equation}
and 
\begin{equation}
\F^{\xi}_k := \sigma (\xi_k, \xi_{k-1}, \ldots) , 
\quad k \in -\N 
\end{equation}
respectively. 
It is obvious that 
\begin{equation}
\F^{\xi}_k \subset \F^{\eta}_k , 
\quad k \in -\N. 
\end{equation}

Since the law on $ G^{-\N} $ of the process $ \eta $ 
determines that of the noise $ \xi $, 
the following definition is reasonable. 

\begin{defi} \label{Def sol}
Let $ \mu = (\mu_k , k \in -\N) $ be a sequence 
of Borel probability measures $ \mu_k $ on $ G $. 
A {\em solution} of the STYE on $ G $ 
with the noise law $ \mu $ 
is a probability measure $ \bP $ on $ G^{-\N} $ such that 
the following two statements hold: 
\\
(i) $ \xi_k $ is independent of $ \F^{\eta}_{k-1} $ under $ \bP $, for any $ k \in -\N $. 
\\
(ii) $ \xi_k $ is distributed as $ \mu_k $ under $ \bP $, for any $ k \in -\N $. 

The totality of solutions of the STYE with the noise law $ \mu $ 
will be denoted by $ \sP_\mu $. 
\end{defi}

\begin{rem} \label{SE remark}
For a given $ \mu $, 
a probability measure $ \bP $ on $ G^{-\N} $ belongs to $ \sP_\mu $ 
if and only if the following (inhomogeneous) Markov property holds: 
\begin{equation}
\bE \ebra{ \phi(\eta_k) \mid \F_{k-1} } 
= \int_{G} \phi( \eta_{k-1} \cdot g ) \mu_k(dg) , \quad k \in -\N 
\label{Markov}
\end{equation}
for any bounded measurable function $ \phi $ on $ G $. 
\end{rem}

In order to give a precise meaning to the table \eqref{table}, 
we need to introduce the notions of strong solution and uniqueness in law. 
We follow the usual terminology in the theory of stochastic differential equations. 

\begin{defi} \label{Def s sol}
Let $ \mu $ be given. 
A solution $ \bP \in \sP_\mu $ is called {\em strong} if 
\begin{equation}
\F^{\eta}_k \subset \F^{\xi}_k \ \mbox{up to $ \bP $-null sets}, \quad k \in -\N . 
\end{equation}
If $ \bP \in \sP_\mu $ is not strong, then it is called {\em non-strong}\footnote
{It is sometimes called {\em weak}.}. 
\end{defi}

\begin{defi} \label{Def uni}
Let $ \mu $ be given. 
It is said that {\em uniqueness in law} holds if 
the set $ \sP_\mu $ consists of at least one element. 
\end{defi}

\

Let $ \mu = (\mu_k , k \in -\N) $ be given 
and consider the set $ \sP_\mu $ of solutions 
of the STYE on $ G $ 
with the noise law $ \mu $. 
Since the condition $ \bP \in \sP_\mu $ is equivalent 
to the Markov property \eqref{Markov}, 
the following is obvious. 

\begin{lem} \label{lem: closed}
The set $ \sP_\mu $ is a closed convex subset of 
the linear topological space of signed measures on $ G^{-\N} $ 
equipped with the weak-star topology. 
\end{lem}

We denote the totality of extremal points of $ \sP_\mu $ by $ \ext $. 
Then Lemma \ref{lem: closed} implies that 
$ \ext \subset \sP_\mu $ 
and any solution $ \bQ \in \sP_\mu $ has an integral representation 
\begin{equation}
\bQ(\cdot) = \int_{\ext} \bP(\cdot) \, \gamma(d\bP) 
\end{equation}
for some Borel probability measure $ \gamma $ on $ \ext $. 

The following is also true. 
\begin{lem} \label{lem: tail}
For $ \bP \in \sP_\mu $, 
the solution $ \bP $ is extremal if and only if 
$ \F^{\eta}_{-\infty } $ is $ \bP $-trivial. 
\end{lem}
We omit the proof, 
since the proof of \cite[Theorem 1, 2)]{Yor} still survives for the STYE's.

\section{Existence of a non-strong solution} \label{sec: exist}

Since the group $ G $ is compact, 
we have the normalized Haar measure $ \nu $: 
That is, there exists a unique positive Borel measure $ \nu $ on $ G $ 
such that $ \nu(G)=1 $ 
and such that 
\begin{equation}
\int_G \phi(g h) \nu(dh) = 
\int_G \phi(h g) \nu(dh) = 
\int_G \phi(h) \nu(dh) 
\end{equation}
for any $ g \in G $ and any bounded measurable function $ \phi $ on $ G $. 

\begin{defi}
A random variable $ U $ is said to be {\em uniformly distributed on $ G $} if 
the law of $ U $ is equal to the normalized Haar measure $ \nu $. 
\end{defi}

Now we restate Theorem \ref{Intro exist}. 

\begin{thm} \label{YorTh}
For an arbitrary noise law $ \mu = \cbra{\mu_k} $, 
there exists a unique element $ \bP^*_{\mu} \in \sP_\mu $ 
whose marginal distributions are uniform: 
\begin{equation}
\lambda_k = \nu , 
\quad k \in -\N. 
\label{YorTh marginal}
\end{equation}
Moreover, under $ \bP^*_{\mu} $, 
\begin{equation}
\mbox{each $ \eta_k $ is independent of the noise $ \F^{\xi}=\sigma(\xi_k : \, k \in -\N) $} . 
\label{YorTh indep}
\end{equation}

In particular, the solution $ \bP^*_{\mu} $ is {\em non-strong} 
unless $ G $ is trivial: i.e. $ G = \{ \unit \} $. 
\end{thm}

Theorem \ref{YorTh} immediately implies the following. 
\begin{coro}
Assume that $ G $ is not trivial. 
Then, for any noise law $ \mu $, 
the case {\Cz} in the table \eqref{table} never occurs. 
\end{coro}

In the sequel, 
we always assume that $ G $ is not trivial. 
Now the problem is how to characterize the trichotomy {\Ca}-{\Cc}. 

\

The key fact to the proof of Theorem \ref{YorTh} is the one-to-one and onto correspondence 
of each solution to what we call an entrance law. 

Let $ \mu = (\mu_k , k \in -\N) $ be given. 
\begin{defi}
A family $ \lambda = (\lambda_k , k \in -\N) $ of probability laws $ \lambda_k $ on $ G $ 
is called an {\em entrance law} for the noise law $ \mu $ 
if the following recurrence relation holds: 
\begin{equation}
\lambda_k = \lambda_{k-1} * \mu_k , 
\quad k \in -\N. 
\label{consistency}
\end{equation}
Here $ \mu * \lambda $ for two measures $ \mu $ and $ \lambda $ on $ G $ 
stands for the convolution of $ \mu $ and $ \lambda $: 
\begin{equation}
\int_G \phi(g) \mu * \lambda (dg) = \int_G \int_G \phi(gh) \mu(dg) \lambda(dh) 
\end{equation}
for any bounded measurable function $ \phi $ on $ G $. 
\end{defi}

For any $ \bP \in \sP_\mu $, 
let $ \lambda = (\lambda_k , k \in -\N) $ denote 
the marginal distributions of $ \eta_k $ on $ G $ , 
i.e., 
\begin{equation}
\lambda_k(\cdot) = \bP (\eta_k \in \cdot) , 
\qquad k \in -\N. 
\label{marginal}
\end{equation}
Then the equation \eqref{STYE} implies that \eqref{consistency} holds. 
Conversely, the following holds. 
\begin{lem} \label{ent law}
Let $ \mu $ be a given noise law. 
Let $ \lambda = (\lambda_k , k \in -\N) $ be 
an entrance law for the noise law $ \mu $: 
The family $ (\lambda_k , k \in -\N) $ satisfies 
the consistency condition \eqref{consistency}. 
Then there exists a unique element $ \bP \in \sP_\mu $ such that \eqref{marginal} holds. 
\end{lem}
This is obvious by Kolmogorov's extention theorem, so we omit the proof. 

\begin{proof}[Proof of Theorem \ref{YorTh}]
By Lemma \ref{ent law}, we see that there exists a solution $ \bP^*_\mu $ 
with marginal distributions given by \eqref{YorTh marginal}. 
Let $ \bE^*_\mu $ denote the expectation with respect to $ \bP^*_\mu $. 
Let $ \phi $ be a bounded measurable function on $ G \times G^{\kbra{k-n+1,\ldots,0}} $. 
Noting that the independence of $ \eta_{k-n} $ and $ \sigma (\xi_{k-n+1},\ldots,\xi_0) $, 
we have 
\begin{eqnarray}
&& \bE^*_\mu \ebra{ \phi \cbra{\eta_k; \xi_{k-n+1}, \ldots, \xi_0} } 
\\
&& = \bE^*_\mu \ebra{ \phi \cbra{\eta_{k-n} \cdot \xi_{k-n+1} \cdots \xi_k; \xi_{k-n+1}, \ldots, \xi_0} } 
\\
&& = \bE^*_\mu \ebra{ \phi \cbra{\eta_{k-n}; \xi_{k-n+1}, \ldots, \xi_0} } . 
\end{eqnarray}
This implies \eqref{YorTh indep}. 
\end{proof}

\section{Homogeneous space structure of the set of extremal points} \label{sec: hom sp}

Let $ \mu = (\mu_k , k \in -\N) $ be given. 
For each $ g \in G $, we define a continuous map $ T_g $ on $ \sP_\mu $ by 
\begin{equation}
T_g(\bP)(\cdot) = \bP \cbra{ g \cdot \eta \in \cdot } , 
\quad \bP \in \sP_\mu 
\end{equation}
where we write 
\begin{equation}
g \cdot \eta = \cbra{ g \cdot \eta_k , k \in -\N } , 
\quad g \in G , \, \eta \in G^{-\N}. 
\end{equation}
Then the following is obvious. 
\begin{lem} \label{Tg action}
The family of maps $ (T_g : \, g \in G ) $ 
defines a group action of $ G $ on $ \sP_\mu $ 
with $ \ext $ an invariant subset, 
i.e., 
\begin{equation}
T_g(\sP_\mu) \subset \sP_\mu 
\quad \mbox{and} \quad 
T_g(\ext) \subset \ext 
\quad 
\mbox{for} \ g \in G 
\label{G action 1}
\end{equation}
and 
\begin{equation}
T_g T_{h^{-1}} = T_{gh^{-1}} \quad \mbox{for $ g,h \in G $}. 
\label{G action 2}
\end{equation}
\end{lem}

The following proposition shows that the point $ \bP^*_\mu $, 
which is defined in Theorem \ref{YorTh}, 
can be considered to be the {\em center} of the set of solutions $ \sP_\mu $. 

\begin{prop} \label{center point}
Let $ \mu $ be given. 
Then, for any $ \bP \in \sP_\mu $, it holds that 
\begin{equation}
\bP^*_\mu = \int_G T_g(\bP) \nu(dg) . 
\label{center point 1}
\end{equation}
\end{prop}

\begin{proof}[Proof of Proposition \ref{center point}]
For any bounded measurable function $ \phi $ on $ G $, we have 
\begin{eqnarray}
&& \int_G \phi(h) \int_G T_g(\bP)(\eta_k \in dh) \nu(dg) 
\\
&& = \int_G \bE \ebra{ \phi(g \cdot \eta_k) } \nu(dg) = \int_G \phi(g) \nu(dg) 
\end{eqnarray}
for any $ k \in -\N $. 
This implies that all marginal distributions of the RHS of \eqref{center point 1} 
are uniform on $ G $. 
\end{proof}

Now we are in a position to give the precise statement of Theorem \ref{Intro action}. 

\begin{thm} \label{Gal1}
Let $ \mu $ be given. 
Then the action of $ (T_g: \, g \in G) $ restricted on $ \ext $ is transitive: 
That is, 
if $ \bP^1 $ and $ \bP^2 $ are two solutions in $ \ext $, 
then there exists an element $ g \in G $ such that $ \bP^1 = T_g (\bP^2) $. 
\end{thm}

Let an extremal point $ \bP^o \in \ext $ be fixed. 
We denote by $ H_\mu(\bP^o) $ the isotropic subgroup at $ \bP^o $: 
\begin{equation}
H_\mu(\bP^o) = \{ g \in G : \, T_g(\bP^o) = \bP^o \} . 
\end{equation}
It is easy to see that $ H_\mu(\bP^o) $ is a closed subgroup of $ G $, 
and hence the quotient set $ G/H_\mu(\bP^o) $ is a compact Hausdorff set. 
Thus Theorem \ref{Gal1} implies the following corollary, 
which reveals the homogeneous space structure of the set $ \ext $. 
\begin{coro}
The action of $ (T_g : \, g \in G) $ restricted on $ \ext $ induces 
a homeomorphism: 
\begin{equation}
G/H_\mu(\bP^o) \stackrel{\sim}{\longrightarrow} \ext. 
\label{homeo}
\end{equation}
\end{coro}

The key to the proof of Theorem \ref{Gal1} is the following lemma, 
which we follow the proof of Yamada--Watanabe's theorem \cite[Proposition 1]{YW}. 
See also \cite[pp. 163--166]{IW}. 

We consider a product space $ G^{-\N} \times G^{-\N} \times G^{-\N} $ 
with its coordinate written as $ (\eta^1,\eta^2,\xi) $. 
Define 
\begin{equation}
\F^{\eta^1,\eta^2}_k 
= \sigma \cbra{ \eta^1_j , \eta^2_j : \, j \le k } , 
\quad k \in -\N . 
\end{equation}

\begin{lem} \label{coupling}
Let $ \bP^1 $ and $ \bP^2 $ be two solutions in $ \sP_\mu $. 
Then there exists a probability measure $ \bQ $ on $ G^{-\N} \times G^{-\N} \times G^{-\N} $ 
such that the following statements hold: 
\\
(i) For $ i=1,2 $, the law on $ G^{-\N} \times G^{-\N} $ of $ (\eta^i,\xi) $ under $ \bQ $ 
coincides with that of $ (\eta,\xi) $ under $ \bP^i $ 
where $ \xi $ is defined in \eqref{STYE2}. 
\\
(ii) For any $ k \in -\N $, 
\begin{equation}
\xi_k \mbox{ is independent of } \F^{\eta^1,\eta^2}_{k-1} \mbox{ under $ \bQ $} . 
\end{equation}
\end{lem}

Let $ \bP^1_{\xi}(\cdot) $ and $ \bP^2_{\xi}(\cdot) $ denote the regular conditional probability 
given $ \F^{\xi}=\sigma(\xi_k: \, k \in -\N) $ such that 
\begin{equation}
\bP^i(\eta \in A , \xi \in B) = \int_B \bP^i_{\xi}(A) \mu(d\xi) 
\end{equation}
for arbitrary measurable sets $ A $ and $ B $ of $ G^{-\N} $, $ i=1,2 $. 
Then the desired probability measure $ \bQ $ is obtained as 
\begin{equation}
\bQ(d\eta^1 \, d\eta^2 \, d\xi) = \bP^1_{\xi}(d\eta^1) \, \bP^2_{\xi}(d\eta^2) \, \mu(d\xi) . 
\end{equation}
We can prove the claims (i) and (ii) 
in the same way as in the proof of Yamada--Watanabe's theorem, 
so we omit the proofs of Lemma \ref{coupling}. 

\begin{proof}[Proof of Theorem \ref{Gal1}]
Let $ \bQ $ be the probability measure given in Lemma \ref{coupling}. 
As far as the end of this paragraph, we omit to write ``$ \bQ $-a.s." 
By Lemma \ref{coupling} (i), we have 
\begin{equation}
\eta^i_k = \eta^i_{k-1} \cdot \xi_k , 
\quad k \in -\N , \ i=1,2. 
\end{equation}
Let $ k \in -\N $ and $ n \in \N $ be arbitrary numbers. 
Then 
\begin{equation}
\eta^i_k = \eta^i_{k-n} \cdot \xi_{k-n+1} \cdots \xi_k , 
\quad i=1,2. 
\end{equation}
Thus we have 
\begin{equation}
\cbra{\eta^1_k} \cdot \cbra{\eta^2_k}^{-1} 
= \cbra{\eta^1_{k-n}} \cdot \cbra{\eta^2_{k-n}}^{-1} , 
\quad k \in -\N , \ n \in \N . 
\label{G action3}
\end{equation}
Since the LHS is irrelevant to $ n \in \N $, 
there exists a random variable $ \eps $ 
which is $ \F^{\eta^1,\eta^2}_{-\infty } $-measurable 
such that 
\begin{equation}
\eps = \cbra{\eta^1_k} \cdot \cbra{\eta^2_k}^{-1} , \quad k \in -\N . 
\label{G action4}
\end{equation}

If we denote by $ \bQ \cbra{ \cdot \mid \eps = g } $ 
the regular conditional probability given $ \eps=g $, 
then we have the following disintegration: 
\begin{equation}
\bQ(\cdot) = \int_G \bQ \cbra{ \cdot \mid \eps = g } \bQ( \eps \in dg ) . 
\label{decomp Q}
\end{equation}
Hence we obtain an integral expression of $ \bP^1 $: 
\begin{equation}
\bP^1(\cdot) = \int_G \bQ \cbra{ \eta^1 \in \cdot \mid \eps = g } \bQ( \eps \in dg ) . 
\end{equation}
By definition, we see that the law 
\begin{equation}
\bQ \cbra{ \eta^1 \in \cdot \mid \eps = g } 
\end{equation}
belongs to $ \sP_\mu $ for $ \bQ( \eps \in dg ) $-almost every $ g \in G $. 
By the assumption that $ \bP^1 $ is extremal, 
we obtain 
\begin{equation}
\bP^1(\cdot) = \bQ \cbra{ \eta^1 \in \cdot \mid \eps = g } 
\label{decomp Q2}
\end{equation}
for $ \bQ( \eps \in dg ) $-almost every $ g \in G $. 
We obtain the similar identity for $ \bP^2 $. 
Note that \eqref{G action4} implies 
\begin{equation}
\bQ \cbra{ \eta^1 = g \cdot \eta^2 \mid \eps = g } = 1 
\quad \mbox{for any $ g \in G $} . 
\label{decomp Q3}
\end{equation}
Therefore we conclude that 
\begin{equation}
\bP^1 = T_g \cbra{\bP^2} 
\quad \mbox{for $ \bQ( \eps \in dg ) $-a.e. $ g \in G $}. 
\end{equation}
This completes the proof. 
\end{proof}

\begin{rem}
In the above proof of Theorem \ref{Gal1}, 
take $ g_0 \in G $ such that $ \bP^1 = T_{g_0} \cbra{\bP^2} $. 
Then it is obvious that the measure $ \bQ( g_0^{-1} \eps \in dg ) $ 
is the normalized Haar measure on $ H_\mu(\bP^2) $. 
This fact leads to the following paradox: 
{\em The tail $ \sigma $-field $ \F^{\eta^1,\eta^2}_{-\infty } $ is 
non-trivial under $ \bQ $, 
whereas both $ \F^{\eta^1}_{-\infty } $ and $ \F^{\eta^2}_{-\infty } $ 
are trivial.} 
\end{rem}

\begin{rem}
Consider the STYE on the direct product group $ G \times G $ 
\begin{equation}
(\eta^1_k,\eta^2_k) = (\eta^1_{k-1},\eta^2_{k-1}) \cdot (\xi^1_k,\xi^2_k) , 
\quad k \in -\N 
\end{equation}
with the noise law $ \tilde{\mu}=(\tilde{\mu}_k : k \in -\N) $ given by 
\begin{equation}
\int_{G \times G} \phi(g,h) \tilde{\mu}_k(dg \times dh) = \int_G \phi(g,g) \mu_k(dg) 
\end{equation}
for any bounded measurable function $ \phi $ on $ G \times G $ and for $ k \in -\N $. 
Denote the set of solutions by $ \tilde{\sP}_{\tilde{\mu}} $. 
We denote the marginal laws of $ (\eta^1,\eta^2) $ 
under the measure $ \bQ $ by $ \tilde{\bP} $. 
Then $ \tilde{\bP} $ 
and the regular conditional probabilities $ \tilde{\bP}(\cdot \mid \eps=g) $ 
belong to $ \tilde{\sP}_{\tilde{\mu}} $. 
Note that \eqref{decomp Q} implies that 
\begin{equation}
\tilde{\bP}(\cdot) = \int_G \tilde{\bP}(\cdot \mid \eps=g) \bQ( \eps \in dg ) . 
\end{equation}
In this integral expression of the solution $ \tilde{\bP} $, 
the integrand $ \tilde{\bP}(\cdot \mid \eps=g) $ 
belongs to $ {\rm ex}(\tilde{\sP}_{\tilde{\mu}}) $ 
for $ \bQ( \eps \in dg ) $-almost every $ g \in G $. 
In fact, by \eqref{decomp Q2} and \eqref{decomp Q3}, 
we see that $ \F^{\eta^1,\eta^2}_{-\infty } $ is $ \tilde{\bP}(\cdot \mid \eps=g) $-trivial 
for $ \bQ( \eps \in dg ) $-almost every $ g \in G $. 
\end{rem}

\section{The subgroup $ H_\mu^s $} \label{sec: H mu s}

To begin with, we recall the well-known {\em Peter--Weyl theorem} for compact groups 
(see, e.g., \cite[Chapter 1]{Sugiura} and \cite[Corollary 13]{JS}). 

Let $ G $ be a compact group and let $ \nu $ denote 
the normalized Haar measure on $ G $. 
Let $ \sG $ denote the totality 
of irreducible unitary representations $ \rho $ of $ G $ 
on a finite-dimensional linear space $ V^{\rho} $. 
Then the following holds: {\em 
The family 
\begin{equation}
( \rho_{i,j} : \, 1 \le i,j \le \dim \rho , \ \rho \in \sG ) 
\end{equation}
forms a total family in the space of continuous functions on $ G $}. 
Here $ ( \rho_{i,j} ) $ 
denotes the matrix element of a representation $ \rho \in \sG $. 

\

In what follows 
we consider the STYE \eqref{STYE} on a compact group $ G $ 
with a fixed noise law $ \mu $. 

To characterize the trichotomy, we introduce a subset $ H^s_\mu $ of $ G $ as follows. 
For an extremal point $ \bP^o \in \ext $ we define 
\begin{equation}
\sH^s_\mu(\bP^o) 
= \Bigl\{ \rho \in \sG : \, 
\text{$ \rho(\eta_k) $ is $ \F^{\xi}_k $-m'ble $ \bP^o $-a.s. for $ k \in -\N $} 
\Bigr\}. 
\end{equation}
Here the word ``m'ble" is abbreviated from ``measurable". 
It is clear from Theorem \ref{Gal1} that 
the set $ \sH^s_\mu (\bP^o) $ 
is independent of the choice of $ \bP^o \in \ext $. 
So we simply write $ \sH^s_\mu (\bP^o) $ as $ \sH^s_\mu $, and define 
\begin{equation}
H^s_\mu 
= \{ g \in G : \, \mbox{$ \rho (g) = {\rm id} $ for every $ \rho \in \sH^s_\mu $} \} . 
\label{H mu s}
\end{equation}

Note that we need to know at least one extremal point $ \bP^o $ 
in order to compute $ H_\mu(\bP^o) $ and $ H^s_\mu $. 
Let us introduce two subsets $ H^1_\mu $ and $ H^2_\mu $ of $ G $ 
which can directly be computed from the noise law $ \mu $ as follows. 
For $ \rho \in \sG $, we set\footnote{
Here the integral in RHS of \eqref{Rk} is interpreted 
as the componentwise integral in a fixed matrix representation of $ \rho $.} 
\begin{equation}
R_k = \int_G \rho(g) \mu_k(dg) 
, \quad k \in -\N . 
\label{Rk}
\end{equation}
Then the following two limits exist for any $ k \in -\N $: 
The first one is 
\begin{equation}
r^1_k[\rho] = 
\lim_{n \to \infty } \| R_{k-n} R_{k-n+1} \cdots R_k \| 
\end{equation}
where $ \| \cdot \| $ denotes 
the operator norm of linear operators on the representation space $ V^{\rho} $. 
The second one is 
\begin{equation}
r^2_k[\rho] = 
\lim_{n \to \infty } \left| \det \cbra{R_{k-n} R_{k-n+1} \cdots R_k} \right| . 
\end{equation}
The convergence of the first limit is obvious by $ \| R_j \| \le 1 $ for any $ j \in -\N $. 
That of the second is ensured by $ |\det R_j| \le 1 $ for any $ j \in -\N $, 
which will be assured by Lemma \ref{lem: det ineq}. 
Now we set 
\begin{equation}
\sH_\mu^i := 
\{ \rho \in \sG : \, r^i_k[\rho] > 0 \ \text{for some $ k \in -\N $} \} 
, \quad i=1,2 
\end{equation}
and define 
\begin{equation}
H^i_\mu := 
\{ g \in G : \, \mbox{$ \rho (g) = {\rm id} $ for every $ \rho \in \sH^i_\mu $} \} 
, \quad i=1,2 , 
\label{H mu}
\end{equation}
where the symbol `$ {\rm id} $' stands for the identity on the representation space $ V^\rho $. 

The following hierarchy is fundamental to our analysis. 
\begin{thm} \label{thm: isot group}
{\rm (i)} 
The three subclasses 
$ \sH^1_\mu $, $ \sH^2_\mu $ and $ \sH^s_\mu $ 
satisfy the following inclusions: 
\begin{equation}
\sH^1_\mu \supset \sH^s_\mu \supset \sH^2_\mu . 
\label{subclass inclusion}
\end{equation}
{\rm (ii)} 
The three subsets 
$ H_\mu^1 $, $ H_\mu^2 $ and $ H^s_\mu $ 
are closed normal subgroups of $ G $ 
such that 
\begin{equation}
H_\mu^1 \subset H_\mu (\bP^o) \subset H^s_\mu \subset H_\mu^2 
\label{isot group inclusion}
\end{equation}
for any $ \bP^o \in \ext $. 
If $ G $ is abelian, then the equalities hold: 
\begin{equation}
H_\mu^1 = H_\mu (\bP^o) = H^s_\mu = H_\mu^2 . 
\label{isot group equality}
\end{equation}
\end{thm}

We remark that this result includes the whole statement of Theorem \ref{Intro isot group}. 
The proof of Theorem \ref{thm: isot group} will be given in section \ref{sec: trichotomy}. 

\begin{rem} \label{isot group another pt}
The isotropic subgroup $ H_\mu(\bP) $ 
at another extremal point $ \bP=T_g(\bP^o) \in \ext $ 
is related to $ H_\mu(\bP^o) $ by 
\begin{equation}
H_\mu(\bP) = g H_\mu(\bP^o) g^{-1} . 
\end{equation}
Hence the isotropic subgroup $ H_\mu(\bP^o) $ 
is not necessarily normal, 
while the subgroup $ H^s_\mu $ is always normal. 
\end{rem}

\section{Examples} \label{sec: example}

\begin{exa} \label{ex: Z/2Z}
Consider the STYE on the group $ \Z/2\Z \simeq \kbra{1,-1} $. 
Since the group $ \Z/2\Z $ is abelian, we have the equalities \eqref{isot group equality}. 
Note that the class $ \sG $ consists of only one element $ \rho $ such that 
\begin{equation}
\rho(1) = 1 , \quad \rho(-1) = -1 . 
\end{equation}
For a noise law $ \mu=(\mu_k: k \in -\N) $, we set $ p_k = \mu_k \cbra{ \kbra{1} } $. 
Now set 
\begin{equation}
r_k = \lim_{n \to \infty } \prod_{j=k-n}^k |2p_j-1| , 
\quad k \in -\N . 
\end{equation}
Then Theorem \ref{Intro trichotomy} leads to the following. 

\begin{prop}
The case {\Ca} or {\Cb} occurs 
according to whether the infinite product $ r_k $ vanishes for any $ k \in -\N $ 
or not. 
\end{prop}

This is obvious, so we omit the proof. 
\end{exa}

\

We give a typical example of the STYE on a non-abelian group 
where $ H_\mu(\bP^o) $ is non-normal 
and hence strictly included in $ H^s_\mu $. 

\begin{exa} \label{ex: sym group}
Consider the symmetric group of degree 3: 
\begin{equation}
\fS_3 = \{ \unit,(12),(23),(13),(123),(132) \} . 
\end{equation}
Set 
\begin{eqnarray}
&& H^o = \kbra{ \unit,(12) } , \\
&& H^1 = (13) H^o = \kbra{ (13),(23) } , \\
&& H^2 = (123) H^o = \kbra{ (123),(132) } . 
\end{eqnarray}
Then $ H^o $ is a non-normal subgroup of $ \fS_3 $ such that 
\begin{equation}
\fS_3 / H^o = \kbra{ H^o,H^1,H^2 } . 
\end{equation}
Let $ \mu = (\mu_k : k \in -\N) $ be the sequence of the uniform laws on $ H^o $: 
$ \mu_k = \nu^o $ for any $ k \in -\N $ where 
\begin{equation}
\nu^o(\kbra{\unit}) = \nu^o(\kbra{(12)}) = 1/2 . 
\label{iid law}
\end{equation}

\begin{prop} \label{prop: sym group}
Consider the STYE on $ \fS_3 $ 
with the noise law $ \mu $ given above. 
Then there exists a solution $ \bP^o \in \sP_\mu $ such that the following hold: 
\\
{\rm (i)} Each $ \eta_k $ under $ \bP^o $ is uniformly distributed on $ H^o $. 
\\
{\rm (ii)} The extremal points $ \ext = \kbra{\bP^o,\bP^1,\bP^2} $, 
where 
\begin{equation}
\bP^1:=T_{(13)}(\bP^o) , \quad \bP^2:=T_{(123)}(\bP^o) . 
\label{sym group 1}
\end{equation}
{\rm (iii)} The isotropic subgroup $ H_\mu(\bP^o) = H^o $. Hence 
\begin{equation}
\{ \unit \} = H_\mu^1 \subsetneq H_\mu(\bP^o) = H^o \subsetneq H_\mu^s = H_\mu^2 = \fS_3 . 
\end{equation}
\end{prop}

\begin{proof}
Note that the family $ \mu = (\mu_k : k \in -\N) $ itself forms an entrance law: 
$ \nu^o = \nu^o * \nu^o $. 
Thus there exists a solution $ \bP^o $ such that each marginal distribution 
$ \bP^o(\eta_k \in \cdot) $ for any $ k \in -\N $ coincides with $ \nu^o $. 
Thus we obtain (i). 

Let $ \bP $ be a solution and let $ \lambda = (\lambda_k : k \in -\N) $ be 
the corresponding entrance law. 
Then we have 
\begin{equation}
\lambda_k(\kbra{g}) 
= \frac{1}{2} \lambda_{k-1}(\kbra{g}) + \frac{1}{2} \lambda_{k-1}(\kbra{g (12)}) , 
\quad g \in \fS_3. 
\end{equation}
This implies that there exist $ p_0,p_1,p_2 \ge 0 $ with $ p_0+p_1+p_2=1 $ 
such that 
\begin{eqnarray}
&& \lambda_k(\kbra{\unit}) = 
   \lambda_k(\kbra{(12)}) = p_0/2 , \\
&& \lambda_k(\kbra{(13)}) = 
   \lambda_k(\kbra{(23)}) = p_1/2 , \\
&& \lambda_k(\kbra{(123)}) = 
   \lambda_k(\kbra{(132)}) = p_2/2 
\end{eqnarray}
for any $ k \in -\N $. Therefore we obtain 
\begin{equation}
\bP = p_0 \bP^o + p_1 \bP^1 + p_2 \bP^2 , 
\end{equation}
where $ \bP^1 $ and $ \bP^2 $ are defined in \eqref{sym group 1}. 
Since the measures $ \bP^o $, $ \bP^1 $ and $ \bP^2 $ are mutually singular, 
we obtain (ii). Hence we obtain (iii). 
This completes the proof. 
\end{proof}
\end{exa}

\section{Proof of the characterization theorem of the trichotomy} \label{sec: trichotomy}

First, we prove Theorem \ref{Intro trichotomy}. 
Before proving it, we need the following. 

\begin{lem} \label{van Kampen1}
The set $ \sH = \sH^s_\mu $ 
is a {\em submodule} of $ \sG $, i.e., the following statements hold: 
\\
(0) If $ \rho_1 \in \sH $ and if $ \rho_2 $ is equivalent to $ \rho_1 $, then $ \rho_2 \in \sH $. 
\\
(i) If $ \rho_1,\rho_2 \in \sH $, then $ \rho_1 \otimes \rho_2 \in \sH $. 
\\
(ii) If $ \rho_1,\rho_2 \in \sH $, then $ \rho_1 \oplus \rho_2 \in \sH $. 
\\
(iii) If $ \rho \in \sH $, then $ \bar{\rho} \in \sH $. 
Here $ \bar{\rho} $ denotes the complex conjugate representation. 
\end{lem}
\noindent
This is obvious, so we omit the proof. 

We utilize the following fact. 
\begin{lem}[van Kampen \cite{Kam}] \label{Kampen}
Let $ \sH $ be a submodule of $ \sG $. 
Suppose that 
\begin{equation}
\mbox{$ \rho(g) = {\rm id} $ for every $ \rho \in \sH $ \ $ \Longrightarrow $ \ $ g=e $ }. 
\end{equation}
Then $ \sH = \sG $. 
\end{lem}

This fact plays a key role in the proof of the {\em Tannaka duality} 
in the representation theory of compact groups. 
For the proof of Lemma \ref{Kampen}, 
see, e.g., \cite[Theorem 13.1]{Tann} and \cite[Lemma 17]{JS}. 

Now we proceed to prove Theorem \ref{Intro trichotomy}. 

\begin{proof}[Proof of Theorem \ref{Intro trichotomy}]
The claim (i) is obvious by definition of the isotropic subgroup $ H_\mu(\bP^o) $. 
The claim (iii) follows immediately from (i) and (ii). 
Thus we need only to prove the claim (ii). 

$ 1^{\circ}). $ 
Suppose that the case {\Cb} occurs, 
i.e., that there exists a strong solution $ \bP \in \sP_\mu $. 
Then it holds that $ \F^{\eta}_{-\infty } \subset \F^{\xi}_{-\infty } $ under $ \bP $. 
Since $ (\xi_k : \, k \in -\N) $ is an independent sequence, 
Kolmogorov's 0-1 law holds so that $ \F^{\eta}_{-\infty } $ is $ \bP $-trivial. 
Then Lemma \ref{lem: tail} says that 
the solution $ \bP $ must be an extremal point: $ \bP \in \ext $. 
Since $ \sH^s_\mu(\bP) = \sG $, we obtain $ H^s_\mu = \kbra{\unit} $. 

Theorem \ref{Gal1} says that all the solutions of $ \bP' \in \ext $ 
are obtained by $ \bP' = T_g(\bP) $ for some $ g \in G $. 
Then it is clear that the solution $ \bP' $ is also strong. 
Therefore we obtain the last claim: All extremal solutions are strong 
and the others non-strong. 

$ 2^{\circ}). $ 
Suppose that $ H^s_\mu = \{ \unit \} $. 
Then we see that 
\begin{equation}
\mbox{$ \rho(g) = {\rm id} $ for every $ \rho \in \sH^s_\mu $ 
\ $ \Longrightarrow $ \ $ g = \unit $}. 
\end{equation}
Applying this fact and Lemma \ref{van Kampen1} to Theorem \ref{Kampen}, 
we conclude that 
$ \sH^s_\mu $ coincides with the whole $ \sG $. 
This shows that all the extremal points are strong and that the case {\Cb} occurs. 
\end{proof}

Second, we prove the inclusions \eqref{subclass inclusion} 
in Theorem \ref{thm: isot group}. 
This is an immediate consequence of the following proposition, 
which generalizes Proposition 2 of \cite{Yor}. 

\begin{prop} \label{Gal3}
Let $ \mu $ be given and let $ \bP \in \sP_\mu $. 
\\
{\rm (i)} 
If $ \rho \in \sH^2_\mu $, then 
\begin{equation}
\bE \ebra{ \rho(\eta_k) \mid \F^{\eta}_{-\infty } \vee \F^{\xi}_k } 
= \rho(\eta_k) 
, \quad k \in -\N . 
\label{eq:3}
\end{equation}
If, moreover, $ \bP \in \ext $, 
then $ \rho(\eta_k) $ is $ \F^{\xi}_k $-measurable $ \bP $-a.s. for any $ k \in -\N $. 
\\
{\rm (ii)} 
If $ \rho \notin \sH^1_\mu $, then 
\begin{equation}
\bE \ebra{ \rho(\eta_k) \mid \F^{\eta}_{-\infty } \vee \F^{\xi}_k } 
= O 
, \quad k \in -\N . 
\label{eq:5}
\end{equation}
If, moreover, $ \bP \in \ext $, 
then $ \rho(\eta_k) $ is {\em never} $ \F^{\xi}_k $-measurable $ \bP $-a.s. for any $ k \in -\N $. 
\end{prop}
\begin{proof}
(i) 
To prove the claim, it suffices to show that 
\eqref{eq:3} for arbitrary small $ k \in -\N $. 

Since $ \rho \in \sH^2_\mu $, 
it holds that 
\begin{equation}
\lim_{k \to -\infty } \prod_{j = k }^{k_0} | \det R_j | > 0 
\end{equation}
for arbitrary small $ k_0 \in -\N $. 
Iterating the equation \eqref{STYE}, we have 
\begin{equation}
\rho(\eta_{k_0}) 
= \rho(\eta_{k_0-n}) \Xi_n , \quad n \in \N 
\end{equation}
where 
\begin{equation}
\Xi_n = \rho(\xi_{k_0-n+1}) \rho(\xi_{k_0-n+2}) 
\cdots \rho(\xi_{k_0}) , \quad n \in \N . 
\end{equation}
Since $ \det \bE \ebra{ \Xi_n } \neq 0 $ for $ n \in \N $, 
we can define 
\begin{equation}
\Phi_n = \cbra{ \bE \ebra{ \Xi_n } }^{-1} \Xi_n , \quad n \in \N . 
\end{equation}
Then the sequence $ (\Phi_n : \, n \in \N) $ 
constitutes a matrix-valued bounded $ (\E^{k_0}_{n}) $-martingale\footnote{
We mean that $ \bE \ebra{ \Phi_{n+1} \mid \E^{k_0}_{n} } = \Phi_{n} $.} 
where 
\begin{equation}
\E^{k_0}_{n} = \sigma ( \xi_{k_0}, \xi_{k_0-1}, \ldots, \xi_{k_0-n+1} ) , 
\quad n \in \N . 
\end{equation}
Therefore $ \Phi_n $ converges to an $ \F^{\xi}_{k_0} $-measurable 
$V^{\rho} \otimes V^{\rho} $-valued random element $ \Phi_{\infty } $ almost surely. 
Since 
\begin{equation}
| \det \Phi_n | = \cbra{ \det \bE \ebra{ \Xi_n } }^{-1} , \quad n \in \N , 
\end{equation}
we obtain $ \det \Phi_{\infty } \neq 0 $ almost surely. 
Taking subsequence if necessary, we see that 
\begin{align}
\Psi_n = \rho(\eta_{k_0-n}) \bE \ebra{ \Xi_n } 
= \rho(\eta_{k_0}) (\Phi_n)^{-1} . 
\end{align}
converges to an $\F^{\eta}_{-\infty }$-measurable 
$V^{\rho} \otimes V^{\rho} $-valued random element 
$ \Psi_{\infty } $ almost surely. 
Therefore we conclude that 
$ \rho(\eta_{k_0}) = \Psi_{\infty } \Phi_{\infty }$ 
is $ \F^{\eta}_{-\infty } \vee \F^{\xi}_{k_0} $-measurable. 

(ii) 
We can easily prove the claim 
by imitating the proof of Proposition 2 of \cite{Yor}. 
So we omit the proof. 
\end{proof}

Third, we prove the rest of Theorem \ref{thm: isot group}. 

\begin{proof}[Proof of Theorem \ref{thm: isot group}]
$ 1^{\circ}). $ 
It is obvious by definition 
that $ H^1_\mu $, $ H^2_\mu $ and $ H^s_\mu $ are closed normal subgroups of $ G $. 

$ 2^{\circ}). $ 
Let $ \bP^o \in \ext $ be fixed. 
We seek an equivalent expression of the condition that $ g \in H_\mu(\bP^o) $. 
Note that $ T_g(\bP^o) = \bP^o $ if and only if 
\begin{equation}
\bP^o (g \cdot \eta_k \in \cdot) = \bP^o (\eta_k \in \cdot) 
, \quad k \in -\N , 
\end{equation}
which is equivalent to 
\begin{equation}
\rho(g) \bE^o [ \rho(\eta_k) ] = \bE^o [ \rho(\eta_k) ] 
, \quad k \in -\N , \ \rho \in \sG 
\label{Gal2 eq1}
\end{equation}
by the Peter--Weyl theorem. 

$ 3^{\circ}). $ 
Suppose that $ g \in H^1_\mu $. 
Let $ \rho \notin \sH^1_\mu $ and $ \bP^o \in \ext $. 
Noting that 
\begin{equation}
\bE^o [ \rho(\eta_k) \mid \F^{\xi}_{k-n} ]
= \rho(\eta_{k-n}) 
R_{k-n+1} R_{k-n+1} \cdots R_k , 
\end{equation}
we have 
\begin{equation}
\left\| \bE^o [ \rho(\eta_k) \mid \F^{\xi}_{k-n} ] \right\| 
\le \| R_{k-n+1} R_{k-n+1} \cdots R_k \| . 
\end{equation}
Letting $ n \to \infty $, we have 
\begin{equation}
\bE^o [\rho(\eta_k) \mid \F^{\xi}_{-\infty } ] = O . 
\end{equation}
Since $ \F^{\xi}_{-\infty } $ is $ \bP^o $-trivial, we obtain 
\begin{equation}
\bE^o [ \rho(\eta_k) ] = O , 
\end{equation}
and then we see that \eqref{Gal2 eq1} holds, 
which proves that $ H^1_\mu \subset H_\mu (\bP^o) $ 
for any $ \bP^o \in \ext $. 

$ 4^{\circ}). $ 
Suppose that $ g \in H_\mu(\bP^o) $. 
Let $ \rho \in \sH^s_\mu $. 
Since $ T_g(\bP^o)=\bP^o $, we have the following identity 
between two joint laws on $ G^{-\N} \times G^{-\N} $: 
\begin{equation}
\bP^o \cbra{ (g \cdot \eta, \xi) \in \cdot } = \bP^o \cbra{ (\eta,\xi) \in \cdot } . 
\end{equation}
Thus we have the following identity 
between two regular conditional distributions on $ G^{-\N} $ given $ \F^{\xi} $ 
(cf. Proof of Lemma \ref{coupling}): 
\begin{equation}
\bP^o_{\xi} \cbra{ g \cdot \eta \in \cdot } 
= \bP^o_{\xi} \cbra{ \eta \in \cdot } 
\quad \mu-\mbox{a.s.} 
\end{equation}
Hence we have 
\begin{equation}
\bE^o_{\xi} \ebra{ \rho(g) \rho(\eta) } 
= \bE^o_{\xi} \ebra{ \rho(\eta) } 
\quad \mu-\mbox{a.s.} 
\end{equation}
Since $ \rho(\eta) $ is $ \F^{\xi} $-measurable, 
we obtain $ \rho(g) \rho(\eta) = \rho(\eta) $, $ \mu $-a.s., 
which implies that $ \rho(g)={\rm id} $. 
Therefore we obtain the inclusion $ H_\mu(\bP^o) \subset H^s_\mu $. 

$ 5^{\circ}). $ 
The inclusion $ H^s_\mu \subset H^2_\mu $ follows from 
the inclusion $ \sH^s_\mu \supset \sH^2_\mu $, 
which is assured by Proposition \ref{Gal3}. 

$ 6^{\circ}). $ 
If $ G $ is abelian, then all irreducible representations of $ G $ are one-dimensional, 
so we obtain $ \sH^1_\mu = \sH^2_\mu $, 
which implies \eqref{isot group equality}. 
\end{proof}

Finally, we prove the following lemma, 
which assures the convergence of $ r^2_k[\rho] $. 

\begin{lem} \label{lem: det ineq}
Let $ \mu_1, \ldots, \mu_n $ be probability measures 
on $ \{ x \in \C^n : \, |x| := (|x_1|^2 + \cdots + |x_n|^2)^{1/2} \le 1 \} $ 
and set 
\begin{equation}
u_i = \int x \mu_i(dx), \quad i = 1, 2, \ldots, n . 
\end{equation}
Then 
\begin{equation}
| \det (u_1 \cdots u_n) | \le 1 . 
\end{equation}
\end{lem}

\noindent
{\em Proof.}\footnote{
The authors are informed of this simple proof by Y.~Takahashi.} 
By Hadamard's inequality, we have 
\begin{equation}
| \det ( u_1,\ldots,u_n ) | \le \prod_i |u_i| . 
\end{equation}
By Jensen's inequality, we have 
\begin{equation}
|u_i| \le \int |x| \mu_i(dx) \le 1 , \quad i = 1, 2, \ldots, n . 
\end{equation}
This completes the proof. 
\qed

\section{Tsirelson--Yor equations} \label{sec: TYE}

Let $ S $ be a Polish space and $ G $ a compact group. 
Let $ \theta : S \to G $ be a measurable map 
and let $ \psi : G \times S \to S $ and $\psi^{-1}:G \times S \to S $ 
be two measurable maps such that 
\begin{equation}\label{TsiYor1}
\psi^{-1}( g, \psi(g,s) ) = \psi( g, \psi^{-1}(g,s) ) = s , 
\quad g \in G , \ s \in S . 
\end{equation}
We consider the following stochastic equation in discrete negative time: 
\begin{equation}
\eta_k = \psi \cbra{ \big. \theta (\eta_{k-1}) , \xi_k }, \quad k \in -\N . 
\tag{\ref{TYE}}
\end{equation}
Moreover, we assume the following. 
\begin{assu} \label{UC} 
The mappings $ \psi $ and $ \theta $ commute in the sense that  
\begin{equation}
\theta ( \psi(g,x) ) = g \cdot \theta (x) , \quad g \in G . 
\end{equation}
\end{assu}

\begin{defi}
Let $ S,G,\psi,\psi^{-1} $ and $ \theta $ as above 
and assume that Assumption \ref{UC} is satisfied. 
Then 
the stochastic equation \eqref{TYE} 
is called a {\em Tsirelson--Yor equation}, 
which will be abbreviated by ``TYE". 
\end{defi}

Following the case of STYE's, 
we introduce the notion of a solution as follows. 
Let $ \eta = ( \eta_k , k \in -\N ) $ denote the coordinate mapping process on $ S^{-\N} $ 
and set 
\begin{equation}
\xi_k = \psi^{-1} \cbra{ \big. \theta \cbra{\eta_{k-1}} , \eta_k }, 
\quad k \in -\N. 
\label{TYE xi}
\end{equation}
The filtrations $ (\F^{\eta}_k) $ and $ (\F^{\xi}_k) $ are defined in the same way. 

\begin{defi} \label{Def sol2}
Let $ \mu = (\mu_k , k \in -\N) $ be a sequence 
of Borel probability measures $ \mu_k $ on $ S $. 
A {\em solution} of the TYE \eqref{TYE} with the noise law $ \mu $ 
is a probability measure $ \bP $ on $ S^{-\N} $ such that 
the following two statements hold: 
\\
(i) $ \xi_k $ is independent of $ \F^{\eta}_{k-1} $ under $ \bP $, for any $ k \in -\N $. 
\\
(ii) $ \xi_k $ is distributed as $ \mu_k $ under $ \bP $, for any $ k \in -\N $. 

The totality of solutions of the TYE \eqref{TYE} with the noise law $ \mu $ 
will be denoted by $ \sP_\mu $. 
\end{defi}

We adopt the same notions of strong solutions and uniqueness in law 
as are defined in Definitions \ref{Def s sol} and \ref{Def uni}. 

If $ S=G $, $ \psi(g,s) $ is the product 
and $ \theta $ is the identity mapping, 
then the TYE \eqref{TYE} is exactly the STYE on $ G $ 
and all the notions of a solution, a strong solution and uniqueness in law 
coincide. 

\

Consider a TYE \eqref{TYE} with a given noise law $ \mu $. 
Denote $ \hat{\eta}_k = \theta(\eta_k) $, 
$ \hat{\mu} = \mu \circ \theta^{-\N} $ and so on. 
Then Assumption \ref{UC} implies that 
\begin{equation}
\hat{\eta}_k = \hat{\eta}_{k-1} \cdot \hat{\xi}_k , 
\quad k \in -\N. 
\label{hat eta}
\end{equation}
This is nothing but the STYE on $ G $ with the noise law $ \hat{\mu} $. 
We define 
\begin{equation}
\hat{\bP} = \text{the law of $ \hat{\eta} $ on $ G^{-\N} $ under $ \bP $}. 
\end{equation}
Then we obtain a mapping 
\begin{equation}
\sP_\mu \ni \bP \mapsto \hat{\bP} \in \hat{\sP}_{\hat{\mu}} 
\label{hat}
\end{equation}
Here $ \hat{\sP}_{\hat{\mu}} $ denotes the set of solutions of 
the STYE \eqref{hat eta} with the noise law $ \hat{\mu} $. 

\

Let us give two examples. 

\begin{exa} \label{ex: Yor eq}
Let $ S=\R $ and $ G=\T = \R/\Z \simeq [0,1) $. 
For $ g \in [0,1) $ and $ s \in \R $, we set 
$ \psi(g,s) = s+g $ and $ \psi^{-1}(g,s) = s-g $. 
Set $ \theta(s) = \alpha (s) $ for $ s \in \R $. 
Then the TYE \eqref{TYE} coincides with the equation 
\begin{equation}
\eta_k = \alpha (\eta_{k-1}) + \xi_k , 
\quad k \in -\N . 
\tag{\ref{Yor eq}}
\end{equation}
Then the equation for $ \hat{\eta_k} = \alpha (\eta_k) $ 
is the STYE on $ \T $, 
which is actually \eqref{STYE torus}. 
\end{exa}

\begin{exa} \label{ex: Tanaka eq}
Let $ S=\R \setminus \kbra{0} $ and $ G = \Z/2\Z \simeq \kbra{1,-1} $. 
For $ g=\pm 1 $ and $ s \in \R \setminus \kbra{0} $, we set 
$ \psi(\pm 1,s) = \pm s $ and $ \theta(s) = \sgn (s) $. 
Then the TYE \eqref{TYE} coincides with the equation 
\begin{equation}
\eta_k = \sgn \cbra{\eta_{k-1}} \cdot \xi_k , 
\quad k \in -\N . 
\tag{\ref{Tanaka}}
\end{equation}
Then the equation for $ \hat{\eta_k} = \sgn \cbra{\eta_k} $ 
is the STYE on $ \Z/2\Z $, 
which is actually \eqref{STYE Z/2Z}. 
\end{exa}

\

Let $ g \in G $. 
For the coordinate process $ (\eta_k : \, k \in -\N) $ on $ S^{-\N} $ 
and the process $ (\xi_k : \, \in -\N ) $ defined by \eqref{TYE xi}, 
we set 
\begin{equation}
\eta'_k := \psi ( g \cdot \theta (\eta_{k-1}) , \xi_k ) , \quad k \in -\N . 
\label{G action1}
\end{equation}
For $ \bP \in \sP_\mu $ for some $ \mu $, 
we define $ T_g(\bP) $ 
by the law of the process $ (\eta'_k : \, k \in -\N) $ under $ \bP $. 

Now we have the following. 

\begin{prop} \label{reduction}
(i) The mapping \eqref{hat} is bijective. 
\\
(ii) $ \bP \in \sP_\mu $ is strong iff so is $ \hat{\bP} \in \hat{\sP}_{\hat{\mu}} $. 
\\
(iii) $ \bP \in \sP_\mu $ is extremal iff so is $ \hat{\bP} \in \hat{\sP}_{\hat{\mu}} $. 
\\
(iv) The family of mappings $ (T_g: \, g \in G) $ defines a group action on $ \sP_\mu $ 
and its restriction on $ \ext \cap \sP_\mu $ is transitive. 
\\
(v) The case {\Ca}, {\Cb} or {\Cc} occurs for the TYE \eqref{TYE} with the noise law $ \mu $ 
iff so does for the STYE on $ G $ with the noise law $ \hat{\mu} $, accordingly. 
\end{prop}

\begin{proof}
Let $ \bP' \in \hat{\sP}_{\hat{\mu}} $ be given. 
For any $ k \in -\N $, 
we define a probability measure $ \pi_k $ on $ (G \times S)^{ \kbra{ k, \ldots, 0 } } $ 
in the following way: 
Let $ (U_k , \xi_k , \xi_{k+1} , \ldots, \xi_0) $ 
be a family of independent random variables 
such that 
$ U_k $ is a $ G $-valued random variable distributed as $ \bP'(\eta_k \in \cdot) $ 
and 
$ \xi_j $ is an $ S $-valued random variable distributed as $ \mu_j $ for $ j=k,\ldots,0 $. 
Set $ \eta_k = \psi ( U_k , \xi_k ) $, 
\begin{equation}
\eta_j = \psi ( \theta (\eta_{j-1}), \xi_j ) , \quad j = k+1,\ldots,0 
\end{equation}
and 
\begin{equation}
U_j = \theta (\eta_j) , \quad j = k+1,\ldots,0 . 
\end{equation}
Then we define 
\begin{equation}
\mbox{$ \pi_k = $ the law of $ ( (U_j,\eta_{j}) : \, j = k, k+1, \ldots, 0 ) $}. 
\end{equation}

Thanks to the consistency assumption \eqref{consistency}, 
we see that the family $ \{ \pi_k : \, k \in -\N \} $ 
satisfies Kolmogorov's consistency condition. 
Therefore Kolmogorov's extension theorem 
ensures the existence of a probability measure $ \bQ $ on $ (G \times S)^{- \N} $ 
whose projection on $ (G \times S)^{ \{ k, \ldots, 0 \} } $ 
coincides with $\pi_k$ for $ k \in -\N $. 
If we define $ \bP $ by the projection of $ \bQ $ on $ S^{-\N} $, 
then we obtain $ \hat{\bP} = \bP' $. 
Therefore we conclude that the mapping \eqref{hat} is surjective. 

The rest of the claims are obvious, so we omit their proofs. 
\end{proof}


\end{document}